\input amstex
\documentstyle{amsppt}
\magnification 1200
\NoRunningHeads
\NoBlackBoxes
\document

\def\J{{\Cal J}}
\def\O{\Cal O}
\def\tg_+{\tilde{\frak g_+}}

\def\Ua{U_q(\tilde\g)}
\def\U2{{\Ua}_2}
\def\g{\frak g}
\def\h{\frak h}

\def\C{\Bbb C}

\def\<{\langle}
\def\>{\rangle}
\def\o{\otimes}
\def\e{\varepsilon}

\def\End{\text{End}}

\topmatter
\title Quantum determinants and quasideterminants
\endtitle
\author {\rm {\bf Pavel Etingof and Vladimir Retakh} \linebreak
\vskip .1in
Department of Mathematics\linebreak
Harvard University\linebreak 
Cambridge, MA 02138, USA\linebreak
and\linebreak
Department of Mathematics\linebreak
University of Arkansas\linebreak
Fayetteville, AR 72701\linebreak
e-mail: etingof\@math.harvard.edu\linebreak
e-mail: vretakh\@comp.uark.edu }
\endauthor
\endtopmatter

\head Introduction\endhead

The notion of a
quasideterminant and a quasiminor of a matrix $A=(a_{ij})$
with not necessarily commuting entries was introduced in \cite{GR1-3}. 
The ordinary determinant of a matrix with commuting entries can be written 
(in many ways) as a product of quasiminors. Furthermore, it was noticed
in \cite {GR1-3, KL, GKLLRT, Mo} that such well-known
noncommutative determinants as the Berezinian, the Capelli
determinant, the quantum determinant of the generating
matrix of the quantum group $U_h(gl_n)$ and the Yangian $Y(gl_n)$ can be
expressed as products of commuting
quasiminors.

The aim of this paper is to extend these results
to a rather general class of 
Hopf algebras given by the Faddeev-Reshetikhin-Takhtajan type relations 
-- the twisted quantum groups defined in Section 1.4. Such quantum groups 
arise when Belavin-Drinfeld classical r-matrices \cite{BD} are quantized. 

Our main result is that the quantum determinant of the 
generating matrix of a twisted quantum group equals 
the product of commuting quasiminors of this matrix.  

{\bf Acknowledgments}
We are indebted to Israel Gelfand for inspiring us to do this work. 
The first author was supported in part by
National Science Foundation and the second author by
Arkansas Science and Technology Authority.

\head 1. Twisted quantum groups\endhead

\subhead 1.1. Quantum $gl_n$\endsubhead

Consider the quantum universal enveloping 
algebra $U=U_h(gl_n)$ \cite{Dr}. This is the h-adically complete
topological Hopf algebra over $\C[[h]]$ 
generated by $E_i,F_i$, $i=1,...,n-1$, and $H_i$, $i=1,...,n$ with 
defining relations 
$$
\gather
[H_i,E_i]=E_i,[H_i,F_i]=-F_i,[H_{i+1},E_i]=-E_i,[H_{i+1},F_i]=F_i,\\
[H_i,E_j]=[H_i,F_j]=0\text{ if } i-j\ne 0,1;[H_i,H_j]=0,\\
[E_i,F_j]=\delta_{ij}\frac{e^{h(H_i-H_{i+1})}-e^{-h(H_i-H_{i+1})}}
{e^h-e^{-h}},\\
E_i^2E_{i\pm1}-(e^h+e^{-h})E_iE_{i\pm 1}E_i+E_{i\pm 1}E_i^2=0,\\
F_i^2F_{i\pm1}-(e^h+e^{-h})F_iF_{i\pm 1}F_i+F_{i\pm 1}F_i^2=0,\\
\tag 1.1\endgather
$$ 
The coproduct, counit, and antipode are defined by 
$$
\gather
\Delta(E_i)=E_i\o e^{h(H_i-H_{i+1})}+1\o E_i,\Delta(F_i)=F_i\o 1+
e^{-h(H_i-H_{i+1})}\o F_i,\\
\Delta(H_i)=H_i\o 1+1\o H_i\\
\e(F_i)=\e(E_i)=\e(H_i)=0, \\ S(E_i)=-E_ie^{-h(H_i-H_{i+1})},
S(F_i)=-e^{h(H_i-H_{i+1})}F_i,S(H_i)=-H_i.\tag 1.2\endgather
$$
Let $U_{\ge 0}$ be the subalgebra of $U$ generated by $E_i$ and $H_i$, 
$U_{\le 0}$ be the subalgebra of $U$ generated by $F_i$ and $H_i$, 
and $U_0$ be the subalgebra generated by $H_i$. 
They are Hopf subalgebras of $U$. Let 
$I_+,I_-$ 
be the kernels of the natural maps $U_{\ge 0}\to U_0$, 
$U_{\le 0}\to U_0$. We also denote by $U_+$ and 
$U_-$ the subalgebras of $U_{\ge 0}$ and $U_{\le 0}$ generated by 
$E_i$ and $F_i$ respectively. 

The Hopf algebra $U$ is quasitriangular:
it admits the universal R-matrix
$$
\Cal R=e^{h\sum_iH_i\o H_i}(1+\sum_{j\ge 1} a_j\o a^j)\in U_{\ge 0}\o 
U_{\le 0},\tag 1.3
$$
where $a_j\in U_+$ and $a^j\in U_-$, and $\e(a_j)=\e(a^j)=0$. 

\subhead 1.2. Twists\endsubhead

\proclaim{Definition 1.1} (Drinfeld) We say that an element $\J \in U\o U$ 
is a twist if $\J =1$ mod $h$, and 
$$
(\e\o 1)(\J )=(1\o \e)(\J )=1,\ 
(\Delta\o 1)(\J )\J _{12}=(1\o\Delta)(\J )\J _{23}.
\tag 1.4
$$
\endproclaim

{\bf Remark.} In this definition and below, $\o$ denotes the tensor 
product completed with respect to the $h$-adic topology. 

\proclaim{Definition 1.2}
We say that a twist $\J $ is upper triangular if
$\J=\J^0\J'$, where $\J^0=e^{h\sum_{ij}a_{ij}H_i\o H_j}\in U_0\o U_0$  
($a_{ij}\in \C[[h]]$), 
and $\J'\in 1+I_+\o I_-$.
\endproclaim

\subhead 1.3. The twisted coproduct\endsubhead

Given any twist $\J $, we define a new coproduct $\Delta_\J (x):=
\J ^{-1}\Delta(x)\J $
on $U$ (from now on we do not use the coproduct $\Delta$
and therefore denote $\Delta_\J$ simply by $\Delta$). 
This coproduct defines a new Hopf algebra structure on $U$. 
We denote the obtained Hopf algebra by $U_\J $. 

The Hopf algebra $U_\J $ is quasitriangular 
with the universal R-matrix
$$
\Cal R_\J =\J _{21}^{-1}\Cal R\J .\tag 1.5
$$

Since $U_\J $ coincides with $U$ as an algebra, it has the same 
representations. 
Let $V$ be the n-dimensional (vector) representation, with the standard basis
$v_i$ such that $H_iv_j=\delta_{ij}v_j$, and $E_iv_{i+1}=v_i$. 
Let $R_\J :V\o V\to V\o V$ be defined by 
$$
R_\J =\Cal R_\J |_{V\o V}.\tag 1.6
$$

\subhead 1.4. The twisted quantum group\endsubhead

Define the quantum function algebra $A_\J$. 
This is the h-adically 
complete algebra over $\C[[h]]$ which is generated 
by $T,T^{-1}\in Mat_n(\C)\o A_\J $ with 
the Faddeev-Reshetikhin-Takhtajan defining relations 
$$
TT^{-1}=T^{-1}T=1,\ 
R_\J ^{12}T^{13}T^{23}=T^{23}T^{13}R_\J ^{12}.\tag 1.7
$$
This algebra 
is a Hopf algebra with $\Delta(T)=T^{12}T^{13},\e(T)=1,S(T)=T^{-1}$. 
We call it {\it the twisted quantum group}.

Let $(,):A_\J\times U_\J\to \C[[h]]$ be the bilinear form defined 
by the formula $(T,x)=\pi_V(x)$ and the properties 
$(1,x)=\e(x)$, $(ab,x)=(a\o b,\Delta(x))$. It is easy to see that this 
form is well defined and satisfies the equation
$(a,xy)=(\Delta(a),x\o y)$. This implies that $(,)$ defines a 
Hopf algebra homomorphism $\theta: U_\J\to A_\J^*$
and a Hopf algebra homomorphism $\theta':A_\J\to U_\J^*$. 

One can show that the map $\theta'$ is injective. 
This property is proved by considering the quasiclassical limit, 
and will be used in Section 3.  

\subhead 1.5. Quantum determinant\endsubhead

We have $T=\sum E_{ij}\o t_{ij}$, where $E_{ij}$ are elementary matrices 
and $t_{ij}\in A_\J $. So we can think of $T$ as the matrix $(t_{ij})$ 
over $A_\J $. Let us define the quantum determinant of this matrix. 

It is known that the Hopf algebra $A_\J $ is a flat deformation 
of the function algebra $\O(GL_n)$. Moreover, it is isomorphic 
to $\O(GL_n)[[h]]$ as a coalgebra (by a map that equals 1 modulo $h$). 
So right $A_\J $-comodules 
correspond to left $GL_n$-modules. Let $\text{Det}$ 
be the 1-dimensional $A_\J $-comodule 
corresponding to the determinant character of $GL_n$. Let $v$ be a 
generator of $\text{Det}$, and $\pi^*:\text{Det}\to \text{Det}\o A_\J $ 
the coaction. We have $\pi^*(v)=v\o D$. The element $D$ is 
obviously grouplike. 
This element is called the quantum determinant of $T$. 
It equals the ordinary determinant modulo $h$. 

\head 2. Quasideterminants and the main theorem.\endhead

\subhead 2.1. Quasideterminants\endsubhead

Quasideterminants were introduced in \cite{GR1}, as follows. 
Let $X$ be an $m\times m$-matrix over an algebra $A$. 
For any $1\le i,j\le m$, let 
$r_i(X)$, $c_j(X)$ be the i-th row and the j-th column of $X$. 
Let $X^{ij}$ be the submatrix of $X$ obtained by removing 
the i-th row and the j-th column from $X$. For a row 
vector $r$ let $r^{(j)}$ be $r$ without the j-th entry.
For a column vector $c$ let $c^{(i)}$ be $c$ without the i-th entry. 
Assume that $X^{ij}$ is invertible. Then the quasideterminant 
$|X|_{ij}\in A$ is defined by the formula
$$
|X|_{ij}=x_{ij}-r_i(X)^{(j)}(X^{ij})^{-1}c_j(X)^{(i)},\tag 2.1
$$
where $x_{ij}$ is the $ij$-th entry of $X$.

For any $n\times n$-matrix $X=(x_{ij})$ over an algebra $A$
and any permutation $\sigma\in S_n$, denote 
by $\text{det}_\sigma(X)$  
the expression
$$
\text{det}_\sigma(X)=\mu_\sigma(|X|_{nn},
|X^{nn}|_{n-1,n-1},...,|X^{n..i,n..i}|_{i-1,i-1},...,x_{11}),\tag 2.2
$$
where $X^{n...i,n...i}$ is the matrix 
obtained from $X$ by erasing rows and columns with numbers $i,...,n$, and 
$\mu_\sigma(a_1,...,a_n)=a_{\sigma 1}...a_{\sigma n}$.  

It is easy to see that 
if $Y$ is an upper triangular matrix with ones on the diagonal 
and $Z$ a lower triangular 
matrix with ones on the diagonal then $\text{det}_\sigma(ZXY)=
\text{det}_\sigma(X)$. 

\subhead 2.2. The main theorem\endsubhead

The main result of this paper is the following theorem. 

\proclaim{Main Theorem} Let $\J $ be a upper
triangular twist. Then the factors in (2.2) commute with each other, 
and for any $\sigma\in S_n$  
$$
D=\text{det}_\sigma(T).\tag 2.3
$$
\endproclaim

This theorem for $\J =1$ was formulated in 
\cite{GR1-3} (see Theorem 4.2 in \cite{GR3}) 
and proved in \cite{KL} (Theorem 3.1). 

Concrete examples of upper triangular twists are contained 
in \cite{Ho}. A review of these examples can be found in Section 4.

\head 3. Proof of the main theorem\endhead

First of all, the first statement of the theorem 
(commutativity of the factors) follows from the second one, 
so it is enough to prove the second statement (formula (2.3)). 

Let $L_\J^+=(\pi_V\o 1)(\Cal R_\J)$, where 
$\pi_V:U_\J\to \End(V)$ defines the vector representation of $U_\J$. 
Let $L_\J^-=(\pi_V\o 1)(\Cal R_{\J,21}^{-1})$. If $\J=1$, we will denote 
$L_\J^\pm$ 
simply by $L^\pm$. 

Let $f_\pm:A_\J\to U_\J$ be the algebra homomorphisms defined by 
the formula $f_\pm(T)=L_\J^\pm$. It is easy to check that they are 
well defined and are coalgebra antihomomorphisms. 

Let $f: A_\J\to U_\J\o U_\J$ be the algebra homomorphism defined by 
$$
f(x)=(f_+\o f_-)(\Delta(x)).
$$ 

\proclaim{Proposition 3.1} $f$ is injective. 
\endproclaim

\demo{Proof} It is easy to see that $f_+(x)=(x,\Cal R_\J)$, and 
$f_-(x)=(x,\Cal R_{\J,21}^{-1})$ (this notation means that we take 
the inner product of $x$ 
with the first component of $\Cal R_\J$ and $\Cal R_{\J,21}^{-1}$
and leave the second component intact). 
Therefore, $f(x)=(x,\Cal R_{\J,12}\Cal R_{\J,31}^{-1})$. This implies that 
if $x\in \text{Ker}(f)$ then $(x,y)=0$ for any $y\in U'$, where 
$U'\subset U_\J$ is the saturated subalgebra generated by the left and right 
components of $\Cal R_\J$. 

We claim that $U'=U_\J$. Indeed, $U_\J$ is a quantization 
of the quasitriangular Lie bialgebra $(gl_n,r)$, where 
$r$ is a classical r-matrix on $gl_n$ such that $r^{21}+r=2\sum_{ij}E_{ij}\o 
E_{ji}$. Thus, the components of $r$ generate $gl_n$, i.e. 
$(gl_n,r)$ is a minimal quasitriangular Lie bialgebra. 
This implies that $U_\J$  is a minimal quasitriangular 
Hopf algebra, i.e. $U'=U_\J$. 

To conclude the argument, we recall that the map $\theta'$ is injective. 
This implies 
that if $(x,y)=0$ for all $y\in U_\J$ then $x=0$. Thus, $\text{Ker}(f)=0$, 
as desired.    
$\square$\enddemo

Proposition 3.1 shows that it is enough to prove (2.3) after applying $f$. 
Therefore, the main theorem is a consequence of the following two 
propositions. 

\proclaim{Proposition 3.2}  
$f(D)=Pe^{hH}\o Pe^{-hH}$, where $H=\sum_{i=1}^NH_i$, and 
$P=e^{h\sum_{ij}(a_{ji}-a_{ij})H_i}$. 
\endproclaim

\proclaim{Proposition 3.3} For any $\sigma\in S_n$ one has 
$f(\text{det}_\sigma(T))=
Pe^{hH}\o Pe^{-hH}$. 
\endproclaim

\demo{Proof of Proposition 3.2} 
Since $D$ is 
grouplike, it is enough for us to show that $f_\pm(D)=Pe^{\pm hH}$.

Define a functor $F$ from right $A_\J$-comodules to 
left $U_\J$-modules, as follows. 
Any right $A_\J$-comodule $W$ is also a left $A_\J^*$-module, 
hence the pullback $\theta^*(W)$ 
is a left $U_\J$-module. We set $F(W):=\theta^*(W)$.  

Consider 
the pushforward functors $f_{\pm *}$ from right $A_\J$-modules to left 
$U_\J$-comodules. Consider also the functors $F_\pm$ from 
left $U_\J$-modules to left $U_\J$-comodules given by 
$\pi^*_{F_+(W)}(w)=(\pi_W\o 1)(\Cal R_\J)w^{(1)}$, and 
$\pi^*_{F_-(W)}(w)=(\pi_W\o 1)(\Cal R_{\J,21}^{-1})w^{(1)}$
(here $w^{(1)}$ means $w$ in the first component).
It is easy to see that $F_\pm\circ F=f_{\pm *}$. 

Let $\chi:U\to \C[[h]]$ be the character defined by $\chi(E_i)=\chi(F_i)=0, 
\chi(H_i)=1$ (the determinant character). It is easy to see 
that $F(\text{Det})=\chi$. Indeed, if $\tilde V$ is the standard comodule
over $A_\J$, then $F(\tilde V)=V$, and $\text{Det},\chi$ are the unique 
1-dimensional subobjects in $\tilde V^{\o n}$ and $V^{\o n}$, respectively. 

Now, $f_\pm(D)$ is the element of $U_\J$ which corresponds to the 
1-dimensional comodule $f_\pm(\text{Det})=F_\pm(\chi)$. This implies that 
$$
f_+(D)=(\chi\o 1)(\Cal R_\J), f_-(D)=(\chi\o 1)(\Cal R_{\J,21}^{-1}).\tag 3.1
$$  
Now the Proposition follows from formula (1.3).  
$\square$\enddemo  

\demo{Proof of Proposition 3.3} We have 
$$
f(T)=(f_+\o f_-)(T^{12}T^{13})=\pi_V^1(\J_{21}^{-1}\Cal R_{12}\J_{12}
\J_{31}^{-1}\Cal R_{31}^{-1}\J_{13}),\tag 3.2
$$ 
where $\pi_V^1$ is $\pi_V$ evaluated in the first component. 
By (1.4), we have 
$$
\J_{12}\J_{31}^{-1}=\J_{3,12}^{-1}\J_{31,2}.\tag 3.3
$$
(Here $\J_{3,12}$ means that the first component of $\J$ acts in
the third component of the tensor product, and the second component 
of $\J$ acts in the first two components of the tensor product, and 
$\J_{31,2}$ is defined similarly). 
Thus, (3.2) implies
$$
\gather
f(T)=\pi_V^1(\J_{21}^{-1}\Cal R_{12}\J_{3,12}^{-1}
\J_{31,2}\Cal R_{31}^{-1}\J_{13})=\\
\pi_V^1(\J_{21}^{-1}\J_{3,21}^{-1}\Cal R_{12}
\Cal R_{31}^{-1}\J_{13,2}\J_{13}).\tag 3.4\endgather
$$ 

It is easy to see that $(\pi_V\o 1)(\J')$ is an upper triangular matrix 
with ones on the diagonal, and 
$(\pi_V\o 1)((\J_{21}')^{-1})$ is a lower triangular matrix 
with ones on the diagonal. 
 
Taking this into account, we obtain
$$
f(\text{det}_\sigma(T))=
\text{det}_\sigma[\pi_V^1((\J_{21}^0)^{-1}(\J_{3,21}^0)^{-1}
\Cal R_{12}\Cal R_{31}^{-1}\J_{13,2}^0\J_{13}^0)].\tag 3.5
$$ 

Recall that $\J^0=e^{h\sum_{ij}a_{ij}H_i\o H_j}$. 
Substituting this into (3.5), we get 
$$
\gather
f(\text{det}_\sigma(T))=\text{det}_\sigma\biggl[\text{diag}(e^{-h
\sum_ia_{ij}(H_i\o 1+1\o H_i)})
e^{-h\sum_{ij}a_{ij}H_j\o H_i}\times \\ L_{12}^+L_{13}^-
\text{diag}(e^{h\sum_ia_{ji}(H_i\o 1+1\o H_i)})
e^{h\sum_{ij}a_{ij}H_j\o H_i}\biggr].\tag 3.6\endgather
$$
Using the fact that all diagonal quasiminors of $L_{12}^+L_{13}^{-}$ 
are of weight zero, we obtain from (3.6): 
$$
f(\text{det}_\sigma(T))=(P\o P)
e^{-h\sum_{ij}a_{ij}H_i\o H_j}\text{det}_\sigma(L_{12}^+L_{13}^-)
e^{h\sum_{ij}a_{ij}H_i\o H_j}.\tag 3.7
$$
By the Main theorem for $\J=1$, we have 
$\text{det}_\sigma(L_{12}^+L_{13}^-)=e^{hH}\o e^{-hH}$. This implies that 
$f(\text{det}_\sigma(T))=Pe^{hH}\o Pe^{-hH}$, as desired. 
$\square$\enddemo

\head 4. Construction of triangular twists\endhead

In this section we will explain a construction of triangular twists 
following the paper of Hodges \cite{Ho}. 

Let $\Gamma_1$, $\Gamma_2$ be disjoint subsets of 
$\{1,...,n-1\}$, and $\tau:\Gamma_1\to\Gamma_2$ a bijection 
such that $|a-b|=1$ iff $|\tau(a)-\tau(b)|=1$. We denote by 
$U_{\ge 0}^m$ the algebra generated by $H_j$ and $E_i,i\in \Gamma_m$ 
($m=1,2$), and 
by $U_{\le 0}^m$ the algebra generated by $H_j$ and $F_i,i\in \Gamma_m$. 
We also denote by $U^m$ the algebra generated by $H_j$ and 
$E_i,F_i,i\in \Gamma_m$. 

Let $\h$ be the linear span of $H_j$. We have $\h=\h_m\oplus \h_m^\perp$, 
where $\h_m$ is the span of $H_i-H_{i+1}$ for $i\in \Gamma_m$ and 
$\h_m^\perp$ is the orthogonal complement of $\h_m$ with respect 
to the standard inner product. Slightly abusing notation, 
we denote by $\tau$ the linear map 
$\h\to\h$ such that $\tau(H_i-H_{i+1})=H_{\tau(i)}-H_{\tau(i)+1}$, 
for $i\in \Gamma_1$, and $\tau(\h_1^\perp)=0$.  
Let $f_\tau: U^1\to U^2$ be the homomorphism of Hopf algebras 
defined by the formula $f_\tau(E_i)=E_{\tau(i)}$, 
$f_\tau(F_i)=F_{\tau(i)}$, $f_\tau(H_i)=\tau(H_i)$. 

Let $\Bbb R=e^{h\sum_i H_i\o H_i}(1+\sum_{j\ge 1} a_j\o a^j)$ 
be the universal R-matrix of $U_1$. Here as before $a_j\in U_+,a^j\in U_-$, 
and $\e(a_j)=\e(a^j)=0$. 

Let $\Theta\in\h\o \h$ be a tensor. 
Let 
$$
\J=e^{-h\Theta}(f_\tau\o 1)(\Bbb R)\in U_{\ge 0}^2\o U_{\le 0}^1.\tag 4.1
$$

\proclaim{Proposition 4.1} 
Let $Z=\sum_i \tau(b_i)\o b_i$, where $b_i$ is an orthonormal basis 
of $\h_1$.  
Suppose that $\Theta$ satisfies the following conditions:

(i) $(x\o 1,Z-\Theta)=(1\o \tau(x),Z-\Theta)=0$, $x\in \h_1$;

(ii) $(\tau(x)\o 1+1\o x,\Theta)=0$, $x\in \h_1$.

Then 
the element $\J$ is a upper triangular twist. 
\endproclaim

\demo{Proof} The properties $(\e\o 1)(\J)=1$, $(1\o \e)(\J)=1$
and the triangularity 
are obvious, so it suffices to prove the second relation in (1.4).  

Denote $(f_\tau\o 1)(\Bbb R)$ by $\hat\Bbb R$. 
The hexagon relations for the R-matrix 
give \linebreak $(\Delta\o 1)(\Bbb R)=\Bbb R_{13}\Bbb R_{23}$, and
$(1\o\Delta)(\Bbb R)=\Bbb R_{13}\Bbb R_{12}$.
From them we get 
$$
\gather
(\Delta\o 1)(\J)\J_{12}=
e^{-h(\Theta_{13}+\Theta_{23})}
\hat\Bbb R_{13}\hat\Bbb R_{23}e^{-h\Theta_{12}}\hat\Bbb R_{12},\\
(1\o \Delta)(\J)\J_{23}=e^{-h(\Theta_{13}+\Theta_{12})}
\hat\Bbb R_{13}\hat\Bbb R_{12}e^{-h\Theta_{23}}\hat\Bbb R_{23}.\tag 4.2
\endgather
$$
Now, identity (i) implies that $[e^{-h\Theta_{12}}\hat\Bbb R_{12},
e^{-h\Theta_{23}}\hat\Bbb R_{23}]=0$ (here it is also used 
that the sets $\Gamma_1,\Gamma_2$ are disjoint). Therefore, 
the second identity of (1.4) is equivalent to the equation
$$
e^{-h\Theta_{23}}
\hat\Bbb R_{13}e^{h\Theta_{23}}=
e^{-h\Theta_{12}}
\hat\Bbb R_{13}e^{h\Theta_{12}}.
\tag 4.3
$$
The last equation is equivalent to $[\Theta_{12}-\Theta_{23},
\hat\Bbb R_{13}]=0$, 
which is equivalent to identity (ii). The proposition is proved. 
$\square$\enddemo  

\proclaim{Proposition 4.2} Equations (i) and (ii) have a solution. 
\endproclaim

\demo{Proof} Make a change of variable $Y=Z-\Theta$. The obtained equations 
with respect to $Y$ are:

(i) $(x\o 1,Y)=(1\o \tau(x),Y)=0$, $x\in \h_1$;

(ii) $(\tau(x)\o 1+1\o x,Y)=x+\tau(x)$, $x\in \h_1$

(here we use that $(x,y)=(\tau(x),\tau(y)),x,y\in \h_1$).
The set of solutions of equation (i) is the space $\h_1^\perp\o \h_2^\perp$. 
Let $Y$ be any vector in this space. Define operators 
$a:\h_1\to \h_1^\perp$, $b:\h_1\to\h_2^\perp$ defined by 
$a(x)=(1\o x,Y)$, $b(x)=(\tau(x)\o 1,Y)$. Then equation (ii) 
is equivalent to 
$$
a(x)+b(x)=x+\tau(x).\tag 4.4
$$

Now we will use the following easy lemma. 

{\bf Lemma.} Let $a:\h_1\to\h_1^\perp$, $b:\h_1\to \h_2^\perp$ be any linear 
maps. Then the equations $a(x)=(1\o x,Y), b(x)=(\tau(x)\o 1,Y)$ have a 
solution in $\h_1^\perp\o \h_2^\perp$ if and only if 
$$
(a(x),\tau(y))=(b(y),x).\tag 4.5
$$
for any $x,y\in \h_1$.

{\it it Proof of the Lemma.}
Since $\h_1\cap \h_2=0$, the maps $\h_1\to (\h_2^\perp)^*$,
$\h_2\to (\h_1^\perp)^*$ given by $z\to (z,*)$ are injective. 
The Lemma easily follows from this observation. 

The Lemma implies that for proving the Proposition it suffices to show 
that equations (4.4),(4.5) have a solution. Substituting 
(4.4) into (4.5), we get 
$$
(a(x),\tau(y))=(y+\tau(y)-a(y),x)=(y,x)+(\tau(y),x),\tag 4.6
$$
since $(a(y),x)=0$. Thus, it suffices to show that there exists 
$a:\h_1\to\h_1^\perp$ such that $(a(x),\tau(y))=(x,y+\tau(y))$. 
This is obvious, since, as we mentioned, the natural map 
$\h_2\to(\h_1^\perp)^*$ is injective. 
$\square$\enddemo 

Thus, can construct an upper triangular twist $\J$ corresponding to
any triple $(\Gamma_1,\Gamma_2,\tau)$. From one such twist one may obtain
an affine space of twists using the following proposition. 

\proclaim{Proposition 4.3} Let $\h_0$ be the space of all $y\in \h$ 
such that $(y,x)=(y,\tau(x))$, $x\in \h_1$. Let $\beta\in\Lambda^2\h_0$. 
Let $\J$ be the upper triangular twist constructed above. 
Then $\J_\beta=\J e^{h\beta}$ is also an upper triangular twist. 
\endproclaim

\demo{Proof} As before, the only thing that requires a proof is that 
$\J_\beta$ is a twist. This is equivalent to saying that 
$e^{h\beta}$ is a twist for $U_\J$. This follows from the fact that 
elements of $\h_0$ are primitive in $U_\J$, as the twist $\J$ has weight $0$ 
with respect to $\h_0$. 
\enddemo

In conclusion we discuss the connection of the above constructions 
with the Belavin-Drinfeld classification of quasitriangular structures
on a simple Lie algebra \cite{BD}. This classification states
that the quasitriangular structures on a simple Lie algebra 
are labeled by two types of data -- discrete data and continuous data. 
The discrete data is a triple $(\Gamma_1,\Gamma_2,\tau)$, where 
$\Gamma_1,\Gamma_2$ are Dynkin subdiagrams of the Dynkin diagram 
of the Lie algebra (not necessarily connected), and $\tau:\Gamma_1\to\Gamma_2$ 
is a Dynkin diagram isomorphism
such that for any $\alpha\in \Gamma_1$ there exists $k$ such that 
$\tau^k(\alpha)\notin\Gamma_1$. The continuous data is a point of a 
certain affine space hanging over any fixed discrete data.  
The algebras $U_{\J_\beta}$ for various $\J,\beta$ 
constructed above provide quatizations 
of all quasitriangular structures on $gl_n$ corresponding to
triples $(\Gamma_1,\Gamma_2,\tau)$ with $\Gamma_1,\Gamma_2$ being 
disjoint. 

If $\Gamma_1,\Gamma_2$ are not disjoint, the above method of 
constructing a twist does not 
work, since the left and right components of $e^{-h\Theta}\hat\Bbb R$ 
no longer commute. However, by \cite{EK}, any quasitriangular structure 
can be quantized by means of a suitable twist. We expect that such a twist 
can be chosen to be upper triangular. In this case, the Main theorem 
will generalize to twisted quantum groups corresponding to all 
triples $(\Gamma_1,\Gamma_2,\tau)$.

\end